# Construction of primitive Pythagorean triples and Pythagorean triples with a common multiplier using gnomons

*Aleshkevich Natalia V.*

*Peter the Great St. Petersburg Polytechnic University*

**Abstract**

The traditional construction of primitive Pythagorean triples by the formulas of two independent variables does not allow their ordering. The paper shows a new view on the construction of primitive Pythagorean triples. A method for constructing primitive Pythagorean triples, based on the use of gnomons equal in area to the squares of the legs from a primitive Pythagorean triple, using two dependent variables, is proposed. This enables the ability to build an ordered table of primitive Pythagorean triples.

Working with ordered data, it is possible to approach in a new way the solution of problems that include systems of Pythagorean triples, such as the construction of Eulerian parallelepipeds, the construction of a Perfect Cuboid, or the problem of coloring natural numbers in two colors so that no Pythagorean triple is monochrome (one of the problems of Ramsey theory).

**Keywords:** Pythagorean triples, gnomon, arithmetic progression.

# Introduction

The objects of our consideration are the Pythagorean numbers, also called Pythagorean triples – triples $(x, y, a)$ of natural numbers satisfying the Pythagorean equation:

$$x^2 + y^2 = a^2.$$

The Pythagorean theorem is a fundamental geometric statement: in any right triangle, the area of a square built on the hypotenuse is equal to the sum of the areas of squares built on the legs.

The general solution are the following formulas [1]:

$$y = 2mn; x = m^2 - n^2; \quad a = m^2 + n^2.$$

These formulas describe exactly once every Pythagorean triple $(x, y, a)$, satisfying the condition $\text{GCD}(x, y, a) = 1$. This means that all sides of the Pythagorean triangle are expressed by relatively prime numbers. This triple of numbers is called a primitive Pythagorean triple.

In any primitive Pythagorean triple one of the legs is an even number and the other is an odd number. In this case the hypotenuse $a$ is an odd number. Without loss of generality, we will assume that $x$ is odd and $y -$ even. Under these constraints we can get all primitive Pythagorean triples and only them.

The variables $m$ and $n$, forming primitive Pythagorean triples, were obtained from very abstract considerations and are not related to each other; that is, they are independent.

The task was to find a geometric interpretation of generation of primitive Pythagorean triples; and to, based on the received interpretation, determine the order on the set of primitive Pythagorean triples, their properties, and quantitative estimates.

The paper proposes a method for generating primitive Pythagorean triples using gnomons. A gnomon is a figure whose addition to a square leads to the construction of a larger square.

The paper also considers the description of gnomons using arithmetic progressions. And using gnomons, the construction of Pythagorean triples with a common multiplier is considered.

**1. Building primitive Pythagorean triples using gnomons**

We will consider the construction of consecutive squares, starting with the square of one. The formula of such a construction:

$$(n+1)^2 = n^2 + 2n + 1.$$

To the constructed square with side $n$, we can add a figure whose area is equal to twice the value of the side of the square plus 1: $2n + 1$. This figure, called a gnomon $(G)$, builds the original square to a larger square; the side of which will be equal to $n + 1$. The thickness of the gnomon $(T)$ will be equal to 1. By constructing $k$ such consecutive gnomons, we can construct a square with side $n + k$. We can combine consecutive gnomons with a thickness equal to 1 into one common gnomon with a thickness equal to $k$.

$$x^2 + G = a^2.$$

We need to build a gnomon that is equal in area to some square: $G_y = y^2$

Then we come to the equation:

$$x^2 + y^2 = a^2.$$

We show the construction of squares of a primitive Pythagorean triple using a generating square with side $S = 2tl$ [2]. Here $\text{GCD}(t, l) = 1$. We assume that $l$ is odd and $t$ is of any parity.

Without loss of generality, we will build a square with an even side.

We increase the side $S$ of the generating square by $2t^2$ and build a larger square with side $y = S + 2t^2$ (Fig. 1). In this case, we obtain a gnomon $G$ with a thickness $2t^2$, placed on the generating square.

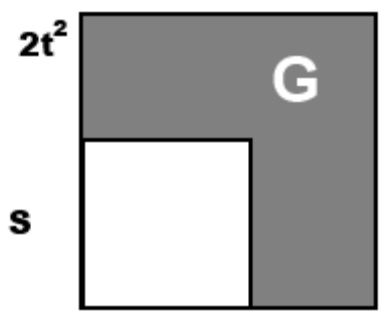

Figure 1

Next, we increase the side $y$ by the value $l^2$. Concurrently we extend the side of the gnomon by the same value $l^2$. At both ends of the gnomon we will have identical rectangles with an area $2t^2 \times l^2$ (Fig. 2).

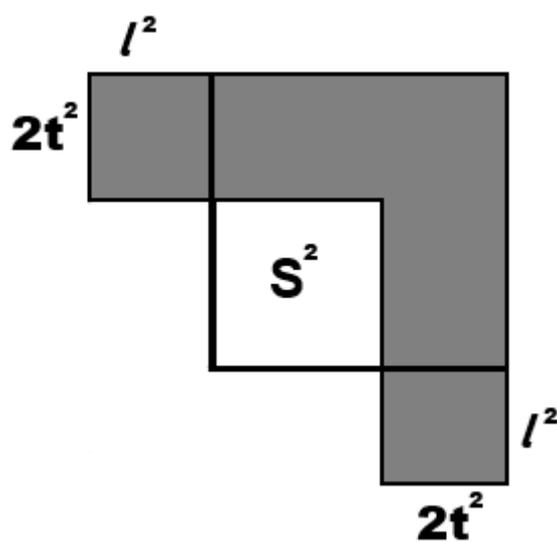

Figure 2

The total area of both rectangles is equal to the area of the generating square with even side. Thus, the newly constructed square with side $y = 2tl + 2t^2 = 2t(t + l)$ unfolds into a gnomon $G_y$ by redistributing the area of the generating square into two equal rectangles, the total area of which is equal to the area of the generating square (Fig. 3). In this case, the area of the gnomon $G_y$ is equal to $y^2$.

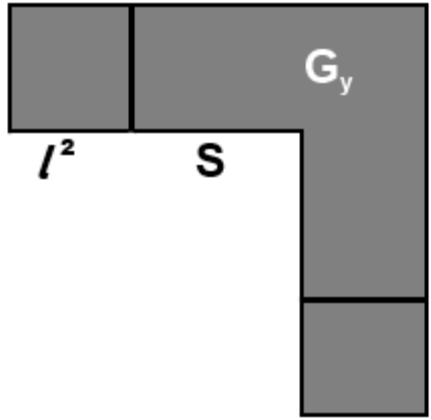

Figure 3

Gnomon $G_y$ is placed on a square with a side:

$$x = S + l^2 = 2tl + l^2 = l(l + 2l) \text{ (Fig. 4)}$$

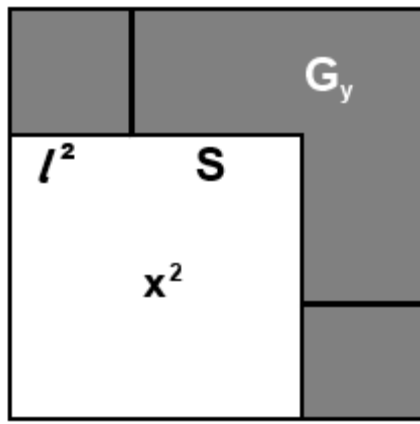

Figure 4

The outer side of the gnomon $G_y$ is equal to the hypotenuse:

$$a = S + 2t^2 + l^2 = (l + t)^2 + t^2.$$

The sum of two squares can be represented as one of the squares and a gnomon placed on it, which is equal in area to the second square. This representation is symmetrical (Fig. 5, 6):

$$x^2 + G_y = y^2 + G_x.$$

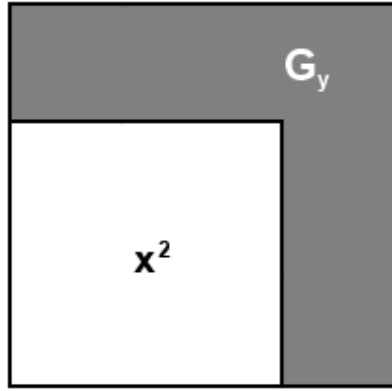

Figure 5

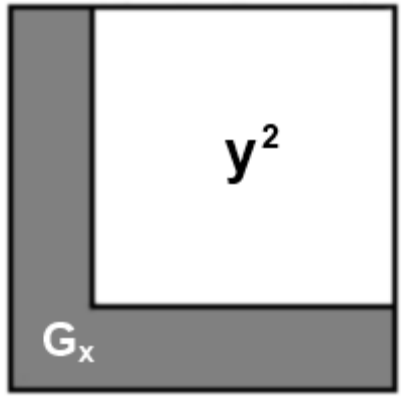

Figure 6

The total square in the form of both gnomons is shown in Fig. 7. Here, the larger gnomon absorbs the smaller gnomon. We will call this representation of gnomons connected gnomons. Both gnomons have a common outer side equals to the hypotenuse $a$. Thus, we have the following relation:

$$x^2 + G_y = y^2 + G_x = G_x + G_y = a^2.$$

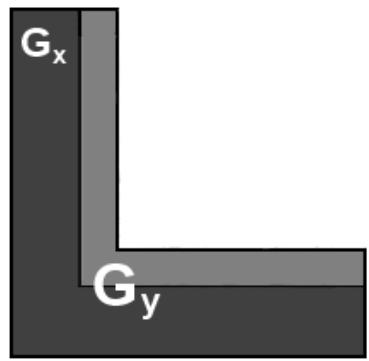

Figure 7

Thus, using generating square, we build gnomons equal in area to the squares of the legs of the primitive Pythagorean triple, which together build the square of the hypotenuse.

**Theorem 1.1. The variables $S, t(S), l(S)$ uniquely determine the primitive Pythagorean triple $(x, y, a)$.**

There are two different proofs of this statement.

*Proof 1*

We represent the area of the generating square by two equal rectangles. For this, we represent the side of the generating square $S$ as a doubled product of two groups of factors:

$$S = 2tl,$$

where

$$t = 2^{\alpha_0 - 1} \prod_{i=1}^{k} p_i^{\alpha_i}$$

$$l = \prod_{i=k+1}^{r} p_i^{\alpha_i},$$

$\alpha_0, \alpha_i$ — powers of the corresponding terms in the product

$p_i$ — various odd prime factors

$$GCD(t, l) = 1.$$

In this case, the number $l$ is always odd, and $t$ can be either an even number when $\alpha_0 > 1$, or an odd number when $\alpha_0 = 1$.

Each of the equal rectangles will have an even side equal to $2t^2$ and an odd side equal to $l^2$. As a result of the construction, we have an even leg from the primitive Pythagorean triple $y = S + 2t^2$, an odd leg from the primitive Pythagorean triple $x = S + l^2$ and a hypotenuse from the primitive Pythagorean triple $a = S + 2t^2 + l^2$.

As we can see, all three numbers of the primitive Pythagorean triple are determined by the values of the variables $S, t(S)$ and $l(S)$. At the same time, these variables are interdependent $S = 2tl$. Together the variables $t$ and $l$ contain all the prime factors into which the value $S$ is factorized. Since every natural number, according to the basic theorem of arithmetic, is factorized into prime factors in the only one way without taking into account the order of the factors, then the obtained numbers defining the Pythagorean triple - $(x, y, a)$ will correspond to the variables $(S, t(S), l(S))$ in the only one way. Theorem 1.1. proved.

*Proof 2*

Expression of the variables $m$ and $n$ through the partitioning variables of the side of the generating square was considered earlier. [2]

We move from our notation of the variables $t$ and $l$ to the generally accepted $m$ and $n$:

$$m = l + t; n = t.$$

In this case, a square with an even side will have a side $y = 2mn$. The total square will have a side $m^2 + n^2$. Indeed, $a = 2tl + 2t^2 + l^2 = t^2 + 2tl + l^2 + t^2 = (t + l)^2 + t^2$.

A square with an odd side will have a side $m^2 - n^2$. Indeed, $x = 2tl + l^2 = 2tl + l^2 + t^2 - t^2 = (t + l)^2 - t^2$.

Since the variables $m$ and $n$, as is known, determine the primitive Pythagorean triple in the only one way, and we use variable substitution to construct similar formulas describing the primitive Pythagorean triple, then using our variables we also describe the primitive Pythagorean triple in the only one way. Theorem 1.1. proved.

**Theorem 1.2. The number of ways to represent the $L(S)$ area of the generating square by two equal rectangles with relatively prime sides, the total area of which is equal to the area of the generating**

square, is determined by the number of odd prime factors without taking into account their powers in the product for the side of the generating square $S = 2 \times 2^{\alpha_0-1} p_1^{\alpha_1} p_2^{\alpha_2} \ldots p_r^{\alpha_r}$, where $p_i$ are prime odd factors, and is equal to

$$L(S) = \sum_{j=0}^{r} C_r^j = 2^r,$$

$r$ — this is the number of prime odd factors without taking into account their powers,

$C_r^j$ — binomial coefficients.

*Proof*

We consider the set of odd prime factors without taking into account their powers of a number $S$: $\{p_1^{\alpha_1}, p_2^{\alpha_2}, \ldots, p_r^{\alpha_r}\}$. Each prime factor is taken into account to an appropriate power as one element. We cannot divide the powers, since the numbers in both factors in the formula must be relatively prime by condition and the factors of these groups represent the sides of a rectangle with relatively prime quantities.

The factors will have partition divided into two groups: $t$ and $l$. In the first step, we do not select any element from the set, then $t$ will be equal to 1, and $l$ will contain all the elements of the set in the product. There is one such sample $C_r^0$.

Next, we can choose any one factor from the list for $t$ and $r - 1$ for $l$, the number of such samples is equal to $C_r^1 = r$. This corresponds to the number of combinations of one element from $r$. Then we can select two elements from the list for $t$ and $r - 2$ elements for $l$. This corresponds to the number of combinations of two elements from : $C_r^2$. And further until the list of odd prime factors is completely ended. In the last step, we select the entire set for $t$ and $l = 1$, this corresponds to the number of combinations of $r$ elements from $r$: $C_r^r$.

After we made all the partitions of the set, we attribute the factor $2^{\alpha_0 - 1}$ at $\alpha_0 > 1$ only to the variable $t$. It is not included in the calculation for binomial coefficients, but is attributed to $t$ for each partition variant using the concatenation operation, as a factor. Thus, the number of partitions of the area of the generating square into two equal rectangles with relatively prime sides with a total area equal to the area of the generating square depends only on the number of odd prime factors $r$ in the composition of the number $S$ without taking into account their powers and is equal to the sum of the binomial coefficients, that is, to number $2^r$.

Theorem 1.2. proved.

Setting the order on a set of primitive Pythagorean triples was considered earlier. [3]

**Theorem 1.3. The variables $S, t(S)$ are setting the order on the set of primitive Pythagorean triples.**

*Proof*

As follows from Theorem 1.1, the variables $S, t(S), l(S)$ uniquely determine the primitive Pythagorean triple $(x, y, a)$. Since the variable $l(S)$ can be expressed via the other two variables:

$$l(S) = \frac{S}{t(S)},$$

For setting the order, we will leave two of them, namely $S$ and $t(S)$.

We build the variable $S$ in ascending order, starting with $S = 2$, and moving in increments of 2:

$$2, 4, 6, 8, \ldots$$

As follows from Theorem 1.2, each value of $S$ corresponds to $L(S) = 2^r$ pairs of numbers $(t(S), l(S))$, where $r$ is the number of odd prime factors in the factorization of the number $S$ into prime factors.

Inside the block by $S$ we will arrange the values of the variable $t(S)$ in ascending order:

$$t_1 < t_2 < \cdots < t_{2^r}.$$

All these numbers differ from each other.

In accordance with the variables $S, t(S)$, ordered tables of primitive Pythagorean triples can be constructed. In this case, ordinal number of the first level $N$ is determined by the equality:

$$N = S/2.$$

The ordinal number of the second level $n$ changes from 1 to $L(S)$ within $S$ in accordance with the growth of $t(S)$.

We write down a set of pairs of elements

$$\{(S_1, t_1(S_1)), \ldots (S_1, t_{2^r}(S_1)), \ldots, (S_k, t_1(S_k)), \quad \ldots, (S_k, t_{2^m}(S_k)), \ldots\}.$$

Since there is a bijection between the elements $(S_i, t_j(S_i))$ and primitive Pythagorean triples (Theorem 1.1), we can write the elements of a primitive Pythagorean triple in a string with a number $(N.n)$. At the same time

$$N = \frac{S}{2}; n(S) = 1, 2, \ldots, L(S).$$

Thus, Theorem 1.3 is proved.

Setting the order on the set of primitive Pythagorean triples makes it possible to construct a table of primitive Pythagorean triples. The table is presented in Appendix 1, where is a fragment of the beginning of the set for the values $S = 2 \div 500$. Accordingly $N = 1, 2, \ldots, 250$.

## 2. Primitive Pythagorean triples and their representation via arithmetic progressions

Primitive Pythagorean triples and their representation via arithmetic progressions we considered earlier. [3]

**Theorem 2.1. Connected gnomons are uniquely described by corresponding arithmetic progressions.**

*Proof*

Without loss of generality, we imagine a primitive Pythagorean triple in the form of a square with side $x$ and a gnomon $G_y$ placed on it. The thickness of the gnomon $T_y$ is equal to the number of consecutive gnomons with a thickness of 1. The area of the first gnomon $s_1$ with a thickness of 1 is equal to $2x + 1$. The area of the next gnomon with a thickness of 1 is equal to

$$s_2 = 2(x + 1) + 1 = 2x + 3.$$

The area of each subsequent gnomon with a thickness of 1 relative to the previous gnomon with a thickness of 1 is equal to

$$s_k = s_{k-1} + 2.$$

The difference of the areas $d$ between consecutive gnomons with a thickness of 1 is equal to two. Therefore, a gnomon describing the square of an even-numbered leg of a primitive Pythagorean triple can be described by three numbers: $(s_1, d, T_y)$. We imagine the number $s_1$, equal to the area of the first gnomon with a thickness of 1, as the first term of the arithmetic progression. The difference of the areas of two consecutive gnomons $d$ with a thickness of 1 will be the difference of the arithmetic progression. The thickness of the gnomon will correspond to the number of terms of the arithmetic progression. There is a one-to-one correspondence between the three numbers describing the gnomon and the three numbers describing the arithmetic progression.

Each gnomon can be described by the corresponding arithmetic progression. Two gnomons (connected gnomons) representing the squares of the legs of a primitive Pythagorean triple, respectively, are described by two arithmetic progressions with different initial terms, a common step $d = 2$, and the corresponding number of terms in each progression:

$$\begin{cases} s_1 = 2x + 1 \\ d = 2 \\ T_y = 2t^2 \end{cases} \Leftrightarrow \begin{cases} s_1 = 2x + 1 \\ d = 2 \\ n = 2t^2 \end{cases}$$

$$\begin{cases} s_1 = 2y + 1 \\ d = 2 \\ T_x = l^2 \end{cases} \Leftrightarrow \begin{cases} s_1 = 2y + 1 \\ d = 2 \\ n = l^2 \end{cases}$$

This proves Theorem 2.1

*Corollary from Theorem 2.1.* The sum of the terms of the arithmetic progression is equal to the area of the corresponding gnomon.

**Theorem 2.2. There is a bijection between variables defining primitive Pythagorean triples $(x, y, a)$ and variables describing arithmetic progressions corresponding to connected gnomons of primitive Pythagorean triples.**

*Proof*

Connected gnomons uniquely define the primitive Pythagorean triple. The gnomon $G_x$ has an area equal to the square of the leg $x$. The gnomon $G_y$ has an area equal to the square of the leg $y$. The outer side of each gnomon is equal to the hypotenuse $a$. There is a correspondence between the hypotenuse and the legs through the thickness of the gnomons describing the squares of these legs:

$$a = x + T_y = y + T_x$$

Imagine a primitive Pythagorean triple in the form of a square and a gnomon placed on it. We take a square with an odd side $x$. Then the area of the gnomon $G_y$ can be represented as the sum of an arithmetic progression with the first term $2x + 1$. Each subsequent term of the progression will be two units larger than the previous one. The number of such terms in arithmetic progression is equal to the thickness of the gnomon

$$T_y = 2t^2.$$

Now we take a square with an even side $y$. Then the area of the gnomon $G_x$ built on it can be represented as the sum of an arithmetic progression with the first term equal to $2y + 1$. Each subsequent term of the progression will be two units larger than the previous one. The number of such terms in arithmetic progression is equal to the thickness of the gnomon

$$T_x = l^2.$$

We consider two connected gnomons. When $y < x$, all the terms of the arithmetic progression describing the gnomon $G_y$ (and their number is $T_y$), will be equal, respectively, to the last terms in the arithmetic progression representing the gnomon $G_x$. And, conversely, for $x < y$, all the terms of the arithmetic progression describing the gnomon $G_x$ (and their number is $T_x$) will be equal, respectively, to the last terms in the arithmetic progression representing the gnomon $G_y$.

This representation in the form of an arithmetic progression of each gnomon fully corresponds to the figure 7 (Fig. 7) of absorption by a gnomon of a larger area of a connected gnomon of a smaller area.

The sum of the terms of the arithmetic progression is equal to the area of the corresponding gnomon (see Corollary from Theorem 2.1). In turn, the gnomon by construction is equal to the square of the corresponding leg. Therefore, the sum of the terms of the arithmetic progression is equal to the square of this leg. Consequently, the sum of the terms of both arithmetic progressions is equal to the square of the hypotenuse.

The middle term (arithmetic mean) $s_x$ of the arithmetic progression describing the gnomon $G_x$ is equal to the sum of the arithmetic progression divided by the number of its terms:

$$S_x = \frac{x^2}{T_x} = \frac{l^2(l+2t)^2}{l^2} = (l+2t)^2.$$

The middle term (arithmetic mean) $s_y$ of the arithmetic progression describing the gnomon $G_y$ is equal to the sum of the arithmetic progression divided by the number of its terms:

$$S_y = \frac{y^2}{T_y} = \frac{4t^2(l+t)^2}{2t^2} = 2(l+t)^2.$$

The first term of the arithmetic progression can be determined through the middle term (arithmetic mean) and the total number of terms.

Indeed, for a gnomon corresponding to the square of an odd leg, the thickness is an odd number; that is, the number of terms of the corresponding arithmetic progression is an odd number. Therefore, the middle term of the arithmetic progression is

$$s = s_1 + d\left(\frac{T_x - 1}{2}\right) = s_1 + T_x - 1,$$

where $\frac{T_x - 1}{2}$ — the number of terms up to the middle term in the arithmetic progression; the progression step is $d = 2$.

For the gnomon corresponding to the square of an even leg, the thickness is an even number, that is, the number of terms of the corresponding arithmetic progression is an even number. Therefore, the middle term of the arithmetic progression (arithmetic mean) is equal to

$$s = s_1 + d\left(\frac{T_y}{2} - \frac{1}{2}\right) = s_1 + T_y - 1,$$

where $\frac{T_y}{2}$ — the number of terms in the first half of the arithmetic progression;

$\frac{d}{2}$ — the half step of the arithmetic progression determines the shift to find the arithmetic mean, since it is an even number and is in the middle between two odd numbers. All the terms of the arithmetic progression are the odd numbers.

From these equations for the middle term we derive the value for the first term of the arithmetic progression:

$$s_1 = s - (T - 1).$$

Since $s_{1y} = 2x + 1$ and $s_{1x} = 2y + 1$, then we can find the values for $x$ and $y$:

$$x = \frac{s_{1y}-1}{2} = \frac{s_y - T_y}{2} = \frac{2(l+t)^2 - 2t^2}{2} = (l+t)^2 - t^2 = l^2 + 2tl = l(l + 2t).$$

$$y = \frac{s_{1x} - 1}{2} = \frac{s_x - T_x}{2} = \frac{(l+2t)^2 - l^2}{2} = \frac{4t^2 + 4tl}{2} = 2t(l + t).$$

Connected gnomons have the same last terms of arithmetic progression. The last term of the arithmetic progression $s_n$ is equal to the sum of the middle term of the arithmetic progression and the corresponding number of terms in this progression minus one:

$$s_n = s_x + (T_x - 1) = s_y + (T_y - 1).$$

The conclusion is similar to obtaining the first term of an arithmetic progression, but the difference in the first case is replaced by the sum in the second case.

Here, the last term is equal to $s_n = 2a - 1$. Hence the equality for hypotenuse $a$ follows:

$$a = \frac{s_n + 1}{2};$$

$$a = \frac{s_x + T_x}{2} = \frac{(l+2t)^2 + l^2}{2} = \frac{2l^2 + 4t^2 + 4lt}{2} = 2lt + 2t^2 + l^2;$$

$$a = \frac{S_y + T_y}{2} = \frac{2(l+t)^2 + 2t^2}{2} = \frac{2l^2 + 4t^2 + 4lt}{2} = 2lt + 2t^2 + l^2.$$

Substitute $S = 2lt$ into the equations, we have:

$$a = S + 2t^2 + l^2 = x + 2t^2 = y + l^2.$$

The last equation corresponds to the value of hypotenuse $a$, constructed using the generating square with side $S$.

Thus, the result obtained is that the arithmetic progressions corresponding to connected gnomons also represent the primitive Pythagorean triple $(y, x, a)$ in the only one way.

Theorem 2.2 is proved.

## 3. Transformation of Gnomons

The transformation of a gnomon is understood as a simultaneous reduction in the number of terms of the arithmetic progression describing it, with an increase in the middle term of the arithmetic progression (arithmetic mean) by the same number of times. Or vice versa. At the same time, the area of the gnomon remains unchanged. The area of the gnomon is equal to the sum of the terms of the arithmetic progression, describing the gnomon. And besides, the area of the gnomon is equal to the square of the corresponding leg of the primitive Pythagorean triple.

With that the sum of the terms of the arithmetic progression is equal to the area of the corresponding leg, either $y^2$ or $x^2$.

The product of the number of terms of the arithmetic progression by the value of its middle term (the arithmetic mean of the values of all terms) is equal to the area of the corresponding square of the leg. The thickness of the gnomon in the arithmetic progression describing it corresponds to the number of terms of this progression.

**Theorem 3.1. The transformation of a gnomon corresponding to square of a leg of a primitive Pythagorean triple through a change in a thickness of a gnomon leads to the disintegration of the original primitive Pythagorean triple and the construction of a new primitive Pythagorean triple**

*Proof*

Consider the formulas for the legs of the primitive Pythagorean triple:

$$y = 2t(l + t); \quad x = l(l + 2t).$$

The number of terms in the arithmetic progression describing the gnomon $G_y$, representing the square of the even leg $y$, is equal to the thickness of the gnomon $T_y = 2t^2$. The middle term of this arithmetic progression is equal $s_y = 2(l + t)^2$. When the thickness of the gnomon is reduced by $k$ times, the value of the middle term increases by $k$ times. At the same time, the number $k$ is a multiplier in $t$. When the value of the middle term is reduced by $m$ times, the thickness of the gnomon increases by $m$ times. In this case, the number $m$ is a multiplier in $(l + t)^2$.

We write the value of $y$ in the form:

$$y = S + 2t^2.$$

If the value of $y$ is unchanged, a change in the value of $t$ immediately leads to a change in the value of $S$. Various primitive Pythagorean triples correspond to different values of pairs $S, t(S)$ (Theorem 1.1).

The number of terms in the arithmetic progression describing the gnomon $G_x$, which represents the square of the odd leg $x$, is equal to the thickness of the gnomon $T_x = l^2$. The middle term of this arithmetic progression is equal to $s_x = (l + 2t)^2$. When the thickness of the gnomon is reduced by $u$ times, the value of the middle term increases by $u$ times. In this case, the number $u$ is a multiplier in $l$. If we increase the thickness of the gnomon by $v$ times, then the value of the middle term is reduced by $v$ times. In this case, the number $v$ is a multiplier in $(l + 2t)^2$.

We write the value of $x$ in the form:

$$x = S + l^2$$

If the value of $x$ is unchanged, a change in the value of $l$ immediately leads to a change in the value of $S$. Various primitive Pythagorean triples correspond to different values of pairs $S$ and $l(S)$ (Theorem 1.1).

Thus, we proved that during the transformation, the gnomon is placed on another square, different from its second leg from the original primitive Pythagorean triple.

According to the formulas: $a = x + 2t^2$; $a = y + l^2$, the side of the total square will also be equal to another value, since when the value of the second term changes, while the first term remains unchanged, the sum also changes.

Thus, during the transformation of the gnomon, the original primitive Pythagorean triple disintegrates and a new primitive Pythagorean triple is formed.

During the transformation, a gnomon can form a number of new primitive Pythagorean triples. This applies to each connected gnomon. That is, this pair of legs from the primitive Pythagorean triple disintegrates when at least one gnomon is transformed. Theorem 3.1 is proved.

*Corollary from Theorem 3.1.* When the gnomon is transformed, the variable $S$ changes, therefore the block $S$ does not contain the same legs, both even and odd.

**Theorem 3.2. The total number of occurrences of the same leg in various primitive Pythagorean triples as a leg depends on the**

**number of n prime factors without taking into account their powers when factorizating this leg into prime factors and is equal to $2^{n-1}$.**

*Proof*

We write an even leg in the form of a product $y = 2t(l + t)$. Let $n$ be the number of prime factors without taking into account their powers in the factorization of the number $y$ into prime factors.

By construction, the numbers $t$ and $l$ are relatively prime, hence the factors $t$ and $(l + t)$ in the formula for $y$ are also relatively prime. Therefore, the factors included in the variables $t$ and $(l + t)$ are different.

Without loss of generality, we can consider half of an even leg: $y/2 = t(l + t)$, since $t$ and $(l + t)$ have different parity, that is, in the prime factorization there is still a factor of 2 in the corresponding power, and it will not be lost.

For an even leg, the number of its occurrences in various primitive Pythagorean triples as a leg depends on the number of different $t$ in the formula.

Changing the variable $t$ for an even leg $y$ leads to the transformation of the gnomon $G_y$, which represents the square of this leg, since this variable determines the thickness of the gnomon.

Changing the thickness of the gnomon leads to the disintegration of the original primitive Pythagorean triple and the formation of a new primitive Pythagorean triple (Theorem 3.1).

If we select some subset from the n-set, then we naturally will have a partition of the n-set into two parts: the first part is formed by the selected elements for $t$, the second by the remaining ones for $(l + t)$.

Each such partition is obtained exactly twice: since at some point we will definitely choose for $t$ those elements that were left for $(l + t)$ before that.

However, we are only interested in those partitions where the even leg has $t < (l + t)$.

Thus, the number of different primitive Pythagorean triples for the same leg $y = t(l + t)$ is determined by the number of different $t$ that can be substituted into this formula for $y$.

The number of such partitions is half of the total number of partitions. When selecting a null subset, the variable $t$ is assigned one, and all elements of the set are written in $(l + t)$.

Since the number of all subsets of a finite set consisting of $n$ elements is $2^n$, the number of subsets for the multiplier $t$ will be $2^{n-1}$.

It is also the same for the legs $x$. We write the odd leg as a product of $x = l(l + 2t)$. Let $n$ be the number of prime factors without taking into account their powers in the factorization of the number $x$ into prime factors.

By construction, the numbers $t$ and $l$ are relatively prime, hence the factors $l$ and $(l + 2t)$ in the formula for $x$ are also relatively prime. Consequently, the factors included in the variables $l$ and $(l + 2t)$ are different.

Changing the variable $l$ for an odd leg $x$ leads to the transformation of the gnomon $G_x$, which represents the square of this leg, since this variable determines the thickness of the gnomon.

Changing the thickness of the gnomon leads to the disintegration of the original primitive Pythagorean triple and the formation of a new primitive Pythagorean triple (Theorem 3.1).

We will determine how many different primitive triples the same odd leg $x$ can form.

For an odd leg, it depends on the number of different partitions of its multipliers into groups $l$ and $(l + 2t)$.

If we select some subset from the n-set, then we will naturally obtain a partition of the n-set into two parts: the first part is formed by the selected elements for $l$, the second by the remaining ones for $(l + 2t)$.

Each such partition is obtained exactly twice: since at some point we will definitely choose for $l$ those elements that were left for $(l + 2t)$ before that. However, we are only interested in those partitions where the odd leg has $l < (l + 2t)$.

Thus, the number of different primitive Pythagorean triples for the same leg $x = l(l + 2t)$ is determined by the number of different $l$ that can be substituted into this formula for $x$.

The number of such partitions, where $l < (l + 2t)$, is half of the total number of partitions. When selecting a null subset, the variable $l$ is assigned one, and all elements of the set are written in $(l + 2t)$.

Since the number of all subsets of a finite set consisting of $n$ elements is $2^n$, the number of subsets for the multiplier $l$ will be $2^{n-1}$.

Thus, Theorem 3.2 is proved.

*Corollary 1 of Theorem 3.2.* The number of different $S$, which corresponds to the number of different primitive Pythagorean triples for the same leg $y = t(l + t)$, is determined by the number of different $t$ that can be substituted into this formula for $y$.

*Corollary 2 of Theorem 3.2.* The number of different $S$, which corresponds to the number of different primitive Pythagorean triples for the same leg $x = l(l + 2t)$, is determined by the number of different $l$ that can be substituted into this formula for $x$.

# 4. Construction of Pythagorean triples with a common multiplier using gnomons

Earlier we considered the construction of Pythagorean triples with a common multiplier. [3]

We multiply all the elements of a primitive Pythagorean triple by an integer coefficient $k$. Thus we will define a Pythagorean triple with a common multiplier: $(kx, ky, ka)$.

We will build a Pythagorean triple with a common multiplier using a square lattice.

***Theorem 4. In a Pythagorean triple with a common multiplier $(kx, ky, ka)$, the thickness of each gnomon increases $k$ times and the middle terms of arithmetic progressions describing the corresponding gnomons increase $k$ times:***

$$T_{kx} = kT_x; \quad T_{ky} = kT_y; \quad s_{kx} = ks_x; \quad s_{ky} = ks_y.$$

*Proof*

We construct a primitive Pythagorean triple $(x, y, a)$. We draw it as a square with side $a$. Inside this square, in the upper right corner, we will place a square with side $y$. Now we add the gnomon $G_x$ to the inner square. The area of the gnomon is equal to the area of the square with side $x$. Thickness of the gnomon is $T_x = l^2$. In this case, $a = y + l^2$.

Next, we place the square $a^2$ in the cells of the square lattice with side $ka$ (Fig. 8).

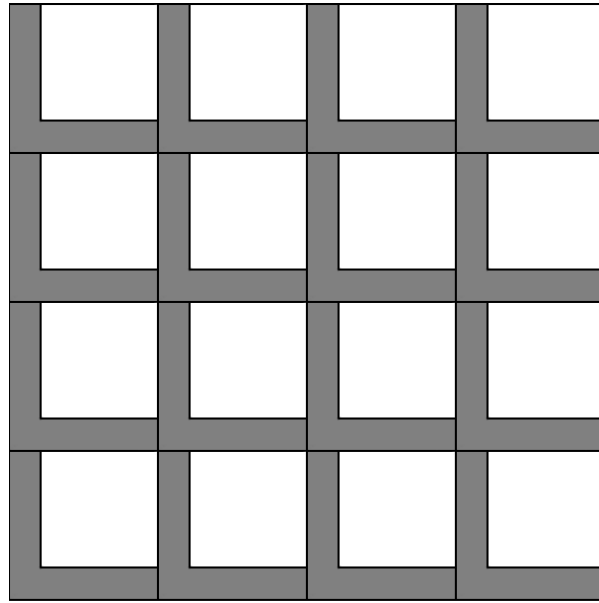
Figure 8

We put together all the squares with side $y$ on the right into a large square with side $ky$, and on the left we will depict the total gnomon assembled from gnomons inside each inner square of the square lattice (Fig. 9). The area of the total square $a^2$ has increased by $k^2$ times. The area of the square with side $y$ has increased by $k^2$ times. The area of the gnomon $G_x$ has also increased by $k^2$ times. Thus, we have constructed a Pythagorean triple with a common multiplier $(kx, ky, ka)$.

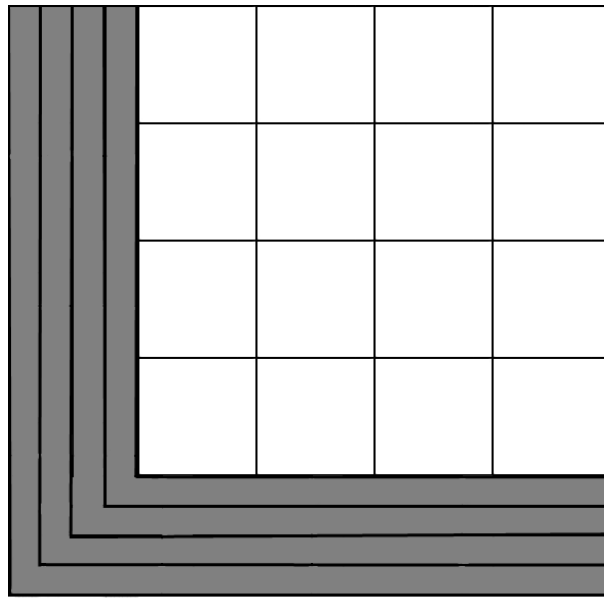
Figure 9

Since the hypotenuse is $ka = k(y + l^2) = ky + kl^2$, then the common side of the collected squares becomes equal to $ky$, and the thickness of the total gnomon $G_{kx}$ becomes equal to $T_{kx} = kl^2 = kT_x$.

The middle term of the arithmetic progression describing the gnomon $G_{kx}$ is calculated by the formula:

$$s_{kx} = \frac{(kx)^2}{T_{kx}} = \frac{k^2 x^2}{kl^2} = k\frac{l^2(l+2t)^2}{l^2} = k(l+2t)^2 = ks_x.$$

The construction will be similar if we place a square with side $x$ in the total square with side $a$ and place a gnomon $G_y$ on it. The thickness of the gnomon is $T_y = 2t^2$.

Next, we will place the square $a^2$ in the cells of the square lattice with side $ka$ (Fig. 8). We will put together all the squares with side $x$ on the right into a large square with side $kx$, and on the left we will depict the total gnomon $G_{ky}$, assembled from gnomons inside each inner square of the square lattice (Fig. 9). The area of the total square $a^2$ has increased by $k^2$ times. The area of the square with side $x$ has increased by $k^2$ times. The area of the gnomon $G_y$ has also increased by $k^2$ times. Thus, we have constructed a Pythagorean triple with a common multiplier $(kx, ky, ka)$.

Since the hypotenuse is $ka = k(x + 2t^2) = kx + 2kt^2$, the common side of the assembled squares becomes equal to $kx$, and the thickness of the total gnomon $G_{ky}$ becomes equal to $T_{ky} = 2kt^2 = kT_y$.

The middle term of the arithmetic progression describing the gnomon $G_{ky}$ is calculated by the formula:

$$s_{ky} = \frac{(ky)^2}{T_{ky}} = \frac{k^2 y^2}{2kt^2} = k\frac{4t^2(l+t)^2}{2t^2} = 2k(l+t)^2 = ks_y.$$

Thus, when multiplying all the elements of the Pythagorean triple by an integer coefficient $k$, the thickness of each gnomon increases by $k$ times and the middle terms of the arithmetic progressions describing the corresponding gnomons increase by $k$ times. Theorem 4 is proved.

**Conclusion**

The construction of primitive Pythagorean triples is based on the concept of a generating square with an even side. With the help of this concept, the relationship of variables necessary for the construction of primitive Pythagorean triples is found. The number of different constructions of primitive Pythagorean triples is determined depending on the size of a side of a generating square.

Abstract formulas for constructing primitive Pythagorean triples through the traditional unrelated parameters $m$ and $n$ are replaced by mutually related parameters through the side of the generating square $S = 2tl$, where $t = n; m = l + t$. A sequential increase in the side of the generating square with a constant step equal to 2, starting from 2, led to the construction of an order on the set of primitive Pythagorean triples using two related parameters $S$ and $t$. A table of primitive Pythagorean triples, constructed in ascending order of the parameter $S$ and the second-level parameter inside the block $S$, namely $t(S)$ up to the value $S = 500$, is given.

Three ways of representing a primitive Pythagorean triple are described: the square of the first leg plus the gnomon of the second, and, conversely, the square of the second leg plus the gnomon of the first leg, as well as two connected gnomons.

The paper shows a description of connected gnomons by arithmetic progressions. Their one-to-one correspondence has been proved. The existence of a bijection between the numbers of a primitive Pythagorean triple and the variables of arithmetic progressions corresponding to connected gnomons is proved.

It has been proven that with the transformation of any gnomon equal in area to the square of the leg, that is, a coordinated change in the thickness of the gnomon and the middle term of the arithmetic progression describing the gnomon, provided that the area of the gnomon is preserved, the leg forms another primitive Pythagorean triple.

A method of assembling gnomons on a square lattice for legs with a common multiplier is proposed. The dependence of the thickness of the gnomon and the middle term of the arithmetic progression describing the gnomon on the value of the common multiplier of the numbers of the Pythagorean triple is proved.

**References**

1. Dunham W. «Euclid's Proof of the Pythagorean Theorem.» Journey Through Genius: The Great Theorems of Mathematics. New York: Wiley Science Editions, (1990) pp. 27-60.

2. Aleshkevich N. V. (2019). Geometric and algebraic interpretation of primitive Pythagorean triples parameters. Retrieved from https://arxiv.org/abs/1907.05540.

3. Aleshkevich N. V. (2021). Constructive representation of primitive Pythagorean triples. Retrieved from https://arxiv.org/abs/2108.06799.

## Fragment of a table of primitive Pythagorean triples constructed with increasing parameters $S, t(S)$

| $N.n_i$ | $S$ | $t$ | $l$ | $x = S + l^2$ | $y = S + 2t^2$ | $a = S + l^2 + 2t^2$ |
|---|---|---|---|---|---|---|
| 1.1 | 2 | 1 | 1 | 3 | 4 | 5 |
| 2.1 | 4 | 2 | 1 | 5 | 12 | 13 |
| 3.1 | 6 | 1 | 3 | 15 | 8 | 17 |
| 3.2 |  | 3 | 1 | 7 | 24 | 25 |
| 4.1 | 8 | 4 | 1 | 9 | 40 | 41 |
| 5.1 | 10 | 1 | 5 | 35 | 12 | 37 |
| 5.2 |  | 5 | 1 | 11 | 60 | 61 |
| 6.1 | 12 | 2 | 3 | 21 | 20 | 29 |
| 6.2 |  | 6 | 1 | 13 | 84 | 85 |
| 7.1 | 14 | 1 | 7 | 63 | 16 | 65 |
| 7.2 |  | 7 | 1 | 15 | 112 | 113 |
| 8.1 | 16 | 8 | 1 | 17 | 144 | 145 |
| 9.1 | 18 | 1 | 9 | 99 | 20 | 101 |
| 9.2 |  | 9 | 1 | 19 | 180 | 181 |
| 10.1 | 20 | 2 | 5 | 45 | 28 | 53 |
| 10.2 |  | 10 | 1 | 21 | 220 | 221 |
| 11.1 | 22 | 1 | 11 | 143 | 24 | 145 |
| 11.2 |  | 11 | 1 | 23 | 264 | 265 |
| 12.1 | 24 | 4 | 3 | 33 | 56 | 65 |
| 12.2 |  | 12 | 1 | 25 | 312 | 313 |
| 13.1 | 26 | 1 | 13 | 195 | 28 | 197 |
| 13.2 |  | 13 | 1 | 27 | 364 | 365 |
| 14.1 | 28 | 2 | 7 | 77 | 36 | 85 |
| 14.2 |  | 14 | 1 | 29 | 420 | 421 |
| 15.1 | 30 | 1 | 15 | 255 | 32 | 257 |
| 15.2 |  | 3 | 5 | 55 | 48 | 73 |
| 15.3 |  | 5 | 3 | 39 | 80 | 89 |
| 15.4 |  | 15 | 1 | 31 | 480 | 481 |
| 16.1 | 32 | 16 | 1 | 33 | 544 | 545 |
| 17.1 | 34 | 1 | 17 | 323 | 36 | 325 |
| 17.2 |  | 17 | 1 | 35 | 612 | 613 |
| 18.1 | 36 | 2 | 9 | 117 | 44 | 125 |
| 18.2 |  | 18 | 1 | 37 | 684 | 685 |
| 19.1 | 38 | 1 | 19 | 399 | 40 | 401 |
| 19.2 |  | 19 | 1 | 39 | 760 | 761 |
| 20.1 | 40 | 4 | 5 | 65 | 72 | 97 |
| 20.2 |  | 20 | 1 | 41 | 840 | 841 |
| 21.1 | 42 | 1 | 21 | 483 | 44 | 485 |
| 21.2 |  | 3 | 7 | 91 | 60 | 109 |
| 21.3 |  | 7 | 3 | 51 | 140 | 149 |
| 21.4 |  | 21 | 1 | 43 | 924 | 925 |
| 22.1 | 44 | 2 | 11 | 165 | 52 | 173 |
| 22.2 |  | 22 | 1 | 45 | 1012 | 1013 |
| 23.1 | 46 | 1 | 23 | 575 | 48 | 577 |
| 23.2 |  | 23 | 1 | 47 | 1104 | 1105 |
| 24.1 | 48 | 8 | 3 | 57 | 176 | 185 |
| 24.2 |  | 24 | 1 | 49 | 1200 | 1201 |
| 25.1 | 50 | 1 | 25 | 675 | 52 | 677 |
| 25.2 |  | 25 | 1 | 51 | 1300 | 1301 |
| 26.1 | 52 | 2 | 13 | 221 | 60 | 229 |

| $N.n_i$ | $S$ | $t$ | $l$ | $x = S + l^2$ | $y = S + 2t^2$ | $a = S + l^2 + 2t^2$ |
|---|---|---|---|---|---|---|
| 26.2 |    | 26 | 1  | 53   | 1404 | 1405 |
| 27.1 | 54 | 1  | 27 | 783  | 56   | 785  |
| 27.2 |    | 27 | 1  | 55   | 1512 | 1513 |
| 28.1 | 56 | 4  | 7  | 105  | 88   | 137  |
| 28.2 |    | 28 | 1  | 57   | 1624 | 1625 |
| 29.1 | 58 | 1  | 29 | 899  | 60   | 901  |
| 29.2 |    | 29 | 1  | 59   | 1740 | 1741 |
| 30.1 | 60 | 2  | 15 | 285  | 68   | 293  |
| 30.2 |    | 6  | 5  | 85   | 132  | 157  |
| 30.3 |    | 10 | 3  | 69   | 260  | 269  |
| 30.4 |    | 30 | 1  | 61   | 1860 | 1861 |
| 31.1 | 62 | 1  | 31 | 1023 | 64   | 1025 |
| 31.2 |    | 31 | 1  | 63   | 1984 | 1985 |
| 32.1 | 64 | 32 | 1  | 65   | 2112 | 2113 |
| 33.1 | 66 | 1  | 33 | 1155 | 68   | 1157 |
| 33.2 |    | 3  | 11 | 187  | 84   | 205  |
| 33.3 |    | 11 | 3  | 75   | 308  | 317  |
| 33.4 |    | 33 | 1  | 67   | 2244 | 2245 |
| 34.1 | 68 | 2  | 17 | 357  | 76   | 365  |
| 34.2 |    | 34 | 1  | 69   | 2380 | 2381 |
| 35.1 | 70 | 1  | 35 | 1295 | 72   | 1297 |
| 35.2 |    | 5  | 7  | 119  | 120  | 169  |
| 35.3 |    | 7  | 5  | 95   | 168  | 193  |
| 35.4 |    | 35 | 1  | 71   | 2520 | 2521 |
| 36.1 | 72 | 4  | 9  | 153  | 104  | 185  |
| 36.2 |    | 36 | 1  | 73   | 2664 | 2665 |
| 37.1 | 74 | 1  | 37 | 1443 | 76   | 1445 |
| 37.2 |    | 37 | 1  | 75   | 2812 | 2813 |
| 38.1 | 76 | 2  | 19 | 437  | 84   | 445  |
| 38.2 |    | 38 | 1  | 77   | 2964 | 2965 |
| 39.1 | 78 | 1  | 39 | 1599 | 80   | 1601 |
| 39.2 |    | 3  | 13 | 247  | 96   | 265  |
| 39.3 |    | 13 | 3  | 87   | 416  | 425  |
| 39.4 |    | 39 | 1  | 79   | 3120 | 3121 |
| 40.1 | 80 | 8  | 5  | 125  | 208  | 233  |
| 40.2 |    | 40 | 1  | 81   | 3280 | 3281 |
| 41.1 | 82 | 1  | 41 | 1763 | 84   | 1765 |
| 41.2 |    | 41 | 1  | 83   | 3444 | 3445 |
| 42.1 | 84 | 2  | 21 | 525  | 92   | 533  |
| 42.2 |    | 6  | 7  | 133  | 156  | 205  |
| 42.3 |    | 14 | 3  | 93   | 476  | 485  |
| 42.4 |    | 42 | 1  | 85   | 3612 | 3613 |
| 43.1 | 86 | 1  | 43 | 1935 | 88   | 1937 |
| 43.2 |    | 43 | 1  | 87   | 3784 | 3785 |
| 44.1 | 88 | 4  | 11 | 209  | 120  | 241  |
| 44.2 |    | 44 | 1  | 89   | 3960 | 3961 |
| 45.1 | 90 | 1  | 45 | 2115 | 92   | 2117 |
| 45.2 |    | 5  | 9  | 171  | 140  | 221  |
| 45.3 |    | 9  | 5  | 115  | 252  | 277  |
| 45.4 |    | 45 | 1  | 91   | 4140 | 4141 |
| 46.1 | 92 | 2  | 23 | 621  | 100  | 629  |
| 46.2 |    | 46 | 1  | 93   | 4324 | 4325 |
| 47.1 | 94 | 1  | 47 | 2303 | 96   | 2305 |
| 47.2 |    | 47 | 1  | 95   | 4512 | 4513 |
| 48.1 | 96 | 16 | 3  | 105  | 608  | 617  |
| 48.2 |    | 48 | 1  | 97   | 4704 | 4705 |

| $N.n_i$ | $S$ | $t$ | $l$ | $x = S + l^2$ | $y = S + 2t^2$ | $a = S + l^2 + 2t^2$ |
|---|---|---|---|---|---|---|
| 49.1 | 98 | 1 | 49 | 2499 | 100 | 2501 |
| 49.2 | | 49 | 1 | 99 | 4900 | 4901 |
| 50.1 | 100 | 2 | 25 | 725 | 108 | 733 |
| 50.2 | | 50 | 1 | 101 | 5100 | 5101 |
| 51.1 | 102 | 1 | 51 | 2703 | 104 | 2705 |
| 51.2 | | 3 | 17 | 391 | 120 | 409 |
| 51.3 | | 17 | 3 | 111 | 680 | 689 |
| 51.4 | | 51 | 1 | 103 | 5304 | 5305 |
| 52.1 | 104 | 4 | 13 | 273 | 136 | 305 |
| 52.2 | | 52 | 1 | 105 | 5512 | 5513 |
| 53.1 | 106 | 1 | 53 | 2915 | 108 | 2917 |
| 53.2 | | 53 | 1 | 107 | 5724 | 5725 |
| 54.1 | 108 | 2 | 27 | 837 | 116 | 845 |
| 54.2 | | 54 | 1 | 109 | 5940 | 5941 |
| 55.1 | 110 | 1 | 55 | 3135 | 112 | 3137 |
| 55.2 | | 5 | 11 | 231 | 160 | 281 |
| 55.3 | | 11 | 5 | 135 | 352 | 377 |
| 55.4 | | 55 | 1 | 111 | 6160 | 6161 |
| 56.1 | 112 | 8 | 7 | 161 | 240 | 289 |
| 56.2 | | 56 | 1 | 113 | 6384 | 6385 |
| 57.1 | 114 | 1 | 57 | 3363 | 116 | 3365 |
| 57.2 | | 3 | 19 | 475 | 132 | 493 |
| 57.3 | | 19 | 3 | 123 | 836 | 845 |
| 57.4 | | 57 | 1 | 115 | 6612 | 6613 |
| 58.1 | 116 | 2 | 29 | 957 | 124 | 965 |
| 58.2 | | 58 | 1 | 117 | 6844 | 6845 |
| 59.1 | 118 | 1 | 59 | 3599 | 120 | 3601 |
| 59.2 | | 59 | 1 | 119 | 7080 | 7081 |
| 60.1 | 120 | 4 | 15 | 345 | 152 | 377 |
| 60.2 | | 12 | 5 | 145 | 408 | 433 |
| 60.3 | | 20 | 3 | 129 | 920 | 929 |
| 60.4 | | 60 | 1 | 121 | 840 | 841 |
| 61.1 | 122 | 1 | 61 | 3843 | 124 | 3845 |
| 61.2 | | 61 | 1 | 123 | 7564 | 7565 |
| 62.1 | 124 | 2 | 31 | 1085 | 132 | 1093 |
| 62.2 | | 62 | 1 | 125 | 7812 | 7813 |
| 63.1 | 126 | 1 | 63 | 4095 | 128 | 4097 |
| 63.2 | | 7 | 9 | 207 | 224 | 305 |
| 63.3 | | 9 | 7 | 175 | 288 | 337 |
| 63.4 | | 63 | 1 | 127 | 8064 | 8065 |
| 64.1 | 128 | 64 | 1 | 129 | 8320 | 8321 |
| 65.1 | 130 | 1 | 65 | 4355 | 132 | 4357 |
| 65.2 | | 5 | 13 | 299 | 180 | 349 |
| 65.3 | | 13 | 5 | 155 | 468 | 493 |
| 65.4 | | 65 | 1 | 131 | 8580 | 8581 |
| 66.1 | 132 | 2 | 33 | 1221 | 140 | 1229 |
| 66.2 | | 6 | 11 | 253 | 204 | 325 |
| 66.3 | | 22 | 3 | 141 | 1100 | 1109 |
| 66.4 | | 66 | 1 | 133 | 8844 | 8845 |
| 67.1 | 134 | 1 | 67 | 4623 | 136 | 4625 |
| 67.2 | | 67 | 1 | 135 | 9112 | 9113 |
| 68.1 | 136 | 4 | 17 | 425 | 168 | 457 |
| 68.2 | | 68 | 1 | 137 | 9384 | 9385 |
| 69.1 | 138 | 1 | 69 | 4899 | 140 | 4901 |
| 69.2 | | 3 | 23 | 667 | 156 | 685 |
| 69.3 | | 23 | 3 | 147 | 1196 | 1205 |

| $N.n_i$ | $S$ | $t$ | $l$ | $x = S + l^2$ | $y = S + 2t^2$ | $a = S + l^2 + 2t^2$ |
|---|---|---|---|---|---|---|
| 69.4 |  | 69 | 1 | 139 | 9660 | 9661 |
| 70.1 | 140 | 2 | 35 | 1365 | 148 | 1373 |
| 70.2 |  | 10 | 7 | 189 | 340 | 389 |
| 70.3 |  | 14 | 5 | 165 | 532 | 557 |
| 70.4 |  | 70 | 1 | 141 | 9940 | 9941 |
| 71.1 | 142 | 1 | 71 | 5183 | 144 | 5185 |
| 71.2 |  | 71 | 1 | 143 | 10224 | 10225 |
| 72.1 | 144 | 8 | 9 | 225 | 272 | 353 |
| 72.2 |  | 72 | 1 | 145 | 10512 | 10513 |
| 73.1 | 146 | 1 | 73 | 5475 | 148 | 5477 |
| 73.2 |  | 73 | 1 | 147 | 10804 | 10805 |
| 74.1 | 148 | 2 | 37 | 1517 | 156 | 1525 |
| 74.2 |  | 74 | 1 | 149 | 11100 | 11101 |
| 75.1 | 150 | 1 | 75 | 5775 | 152 | 5777 |
| 75.2 |  | 3 | 25 | 775 | 168 | 793 |
| 75.3 |  | 25 | 3 | 159 | 1400 | 1409 |
| 75.4 |  | 75 | 1 | 151 | 11400 | 11401 |
| 76.1 | 152 | 4 | 19 | 513 | 184 | 545 |
| 76.2 |  | 76 | 1 | 153 | 11704 | 11705 |
| 77.1 | 154 | 1 | 77 | 6083 | 156 | 6085 |
| 77.2 |  | 7 | 11 | 275 | 252 | 373 |
| 77.3 |  | 11 | 7 | 203 | 396 | 445 |
| 77.4 |  | 77 | 1 | 155 | 12012 | 12013 |
| 78.1 | 156 | 2 | 39 | 1677 | 164 | 1685 |
| 78.2 |  | 6 | 13 | 325 | 228 | 397 |
| 78.3 |  | 26 | 3 | 165 | 1508 | 1517 |
| 78.4 |  | 78 | 1 | 157 | 12324 | 12325 |
| 79.1 | 158 | 1 | 79 | 6399 | 160 | 6401 |
| 79.2 |  | 79 | 1 | 159 | 12640 | 12641 |
| 80.1 | 160 | 16 | 5 | 185 | 672 | 697 |
| 80.2 |  | 80 | 1 | 161 | 12960 | 12961 |
| 81.1 | 162 | 1 | 81 | 6723 | 164 | 6725 |
| 81.2 |  | 81 | 1 | 163 | 13284 | 13285 |
| 82.1 | 164 | 2 | 41 | 1845 | 172 | 1853 |
| 82.2 |  | 82 | 1 | 165 | 13612 | 13613 |
| 83.1 | 166 | 1 | 83 | 7055 | 168 | 7057 |
| 83.2 |  | 83 | 1 | 167 | 13944 | 13945 |
| 84.1 | 168 | 4 | 21 | 609 | 200 | 641 |
| 84.2 |  | 12 | 7 | 217 | 456 | 505 |
| 84.3 |  | 28 | 3 | 177 | 1736 | 1745 |
| 84.4 |  | 84 | 1 | 169 | 14280 | 14281 |
| 85.1 | 170 | 1 | 85 | 7395 | 172 | 7397 |
| 85.2 |  | 5 | 17 | 459 | 220 | 509 |
| 85.3 |  | 17 | 5 | 195 | 748 | 773 |
| 85.4 |  | 85 | 1 | 171 | 14620 | 14621 |
| 86.1 | 172 | 2 | 43 | 2021 | 180 | 2029 |
| 86.2 |  | 86 | 1 | 173 | 14964 | 14965 |
| 87.1 | 174 | 1 | 87 | 7743 | 176 | 7745 |
| 87.2 |  | 3 | 29 | 1015 | 192 | 1033 |
| 87.3 |  | 29 | 3 | 183 | 1856 | 1865 |
| 87.4 |  | 87 | 1 | 175 | 15312 | 15313 |
| 88.1 | 176 | 8 | 11 | 297 | 304 | 425 |
| 88.2 |  | 88 | 1 | 177 | 15664 | 15665 |
| 89.1 | 178 | 1 | 89 | 8099 | 180 | 8101 |
| 89.2 |  | 89 | 1 | 179 | 16020 | 16021 |
| 90.1 | 180 | 2 | 45 | 2205 | 188 | 2213 |

| $N.n_i$ | $S$ | $t$ | $l$ | $x = S + l^2$ | $y = S + 2t^2$ | $a = S + l^2 + 2t^2$ |
|---|---|---|---|---|---|---|
| 90.2 |  | 10 | 9 | 261 | 380 | 461 |
| 90.3 |  | 18 | 5 | 205 | 828 | 853 |
| 90.4 |  | 90 | 1 | 181 | 16200 | 16201 |
| 91.1 | 182 | 1 | 91 | 8463 | 184 | 8465 |
| 91.2 |  | 7 | 13 | 351 | 280 | 449 |
| 91.3 |  | 13 | 7 | 231 | 520 | 569 |
| 91.4 |  | 91 | 1 | 183 | 16744 | 16745 |
| 92.1 | 184 | 4 | 23 | 713 | 216 | 745 |
| 92.2 |  | 92 | 1 | 185 | 17112 | 17113 |
| 93.1 | 186 | 1 | 93 | 8835 | 188 | 8837 |
| 93.2 |  | 3 | 31 | 1147 | 204 | 1165 |
| 93.3 |  | 31 | 3 | 195 | 2108 | 2117 |
| 93.4 |  | 93 | 1 | 187 | 17484 | 17485 |
| 94.1 | 188 | 2 | 47 | 2397 | 196 | 2405 |
| 94.2 |  | 94 | 1 | 189 | 17860 | 17861 |
| 95.1 | 190 | 1 | 95 | 9215 | 192 | 9217 |
| 95.2 |  | 5 | 19 | 551 | 240 | 553 |
| 95.3 |  | 19 | 5 | 215 | 912 | 937 |
| 95.4 |  | 95 | 1 | 191 | 18240 | 18241 |
| 96.1 | 192 | 32 | 3 | 201 | 2240 | 2249 |
| 96.2 |  | 96 | 1 | 193 | 18624 | 18625 |
| 97.1 | 194 | 1 | 97 | 9603 | 196 | 9605 |
| 97.2 |  | 97 | 1 | 195 | 19012 | 19013 |
| 98.1 | 196 | 2 | 49 | 2597 | 204 | 2605 |
| 98.2 |  | 98 | 1 | 197 | 19404 | 19405 |
| 99.1 | 198 | 1 | 99 | 9999 | 200 | 10001 |
| 99.2 |  | 9 | 11 | 319 | 360 | 481 |
| 99.3 |  | 11 | 9 | 279 | 440 | 521 |
| 99.4 |  | 99 | 1 | 199 | 19800 | 19801 |
| 100.1 | 200 | 4 | 25 | 825 | 232 | 857 |
| 100.2 |  | 100 | 1 | 201 | 20200 | 20201 |
| 101.1 | 202 | 1 | 101 | 10403 | 204 | 10405 |
| 101.2 |  | 101 | 1 | 203 | 20604 | 20605 |
| 102.1 | 204 | 2 | 51 | 2805 | 212 | 2813 |
| 102.2 |  | 6 | 17 | 493 | 276 | 565 |
| 102.3 |  | 34 | 3 | 213 | 2516 | 2525 |
| 102.4 |  | 102 | 1 | 205 | 21012 | 21013 |
| 103.1 | 206 | 1 | 103 | 10815 | 208 | 10817 |
| 103.2 |  | 103 | 1 | 207 | 21424 | 21425 |
| 104.1 | 208 | 8 | 13 | 377 | 336 | 505 |
| 104.2 |  | 104 | 1 | 209 | 21840 | 21841 |
| 105.1 | 210 | 1 | 105 | 11235 | 212 | 11237 |
| 105.2 |  | 3 | 35 | 1435 | 228 | 1453 |
| 105.3 |  | 5 | 21 | 651 | 260 | 701 |
| 105.4 |  | 7 | 15 | 435 | 308 | 533 |
| 105.5 |  | 15 | 7 | 259 | 420 | 469 |
| 105.6 |  | 21 | 5 | 235 | 1092 | 1117 |
| 105.7 |  | 35 | 3 | 219 | 2660 | 2669 |
| 105.8 |  | 105 | 1 | 211 | 22260 | 22261 |
| 106.1 | 212 | 2 | 53 | 3021 | 220 | 3029 |
| 106.2 |  | 106 | 1 | 213 | 22684 | 22685 |
| 107.1 | 214 | 1 | 107 | 11663 | 216 | 11665 |
| 107.2 |  | 107 | 1 | 215 | 23112 | 23113 |
| 108.1 | 216 | 4 | 27 | 945 | 248 | 977 |
| 108.2 |  | 108 | 1 | 217 | 23544 | 23545 |
| 109.1 | 218 | 1 | 109 | 12099 | 220 | 12101 |

| $N.n_i$ | $S$ | $t$ | $l$ | $x = S + l^2$ | $y = S + 2t^2$ | $a = S + l^2 + 2t^2$ |
|---|---|---|---|---|---|---|
| 109.2 |  | 109 | 1 | 219 | 23980 | 23981 |
| 110.1 | 220 | 2 | 55 | 3245 | 228 | 3253 |
| 110.2 |  | 10 | 11 | 341 | 420 | 541 |
| 110.3 |  | 22 | 5 | 245 | 1188 | 1213 |
| 110.4 |  | 110 | 1 | 221 | 24420 | 24421 |
| 111.1 | 222 | 1 | 111 | 12543 | 224 | 12545 |
| 111.2 |  | 3 | 37 | 1591 | 240 | 1609 |
| 111.3 |  | 37 | 3 | 231 | 2960 | 2969 |
| 111.4 |  | 111 | 1 | 223 | 24864 | 24865 |
| 112.1 | 224 | 16 | 7 | 273 | 736 | 785 |
| 112.2 |  | 112 | 1 | 225 | 25312 | 25313 |
| 113.1 | 226 | 1 | 113 | 12995 | 228 | 12997 |
| 113.2 |  | 113 | 1 | 227 | 25764 | 25765 |
| 114.1 | 228 | 2 | 57 | 3477 | 236 | 3485 |
| 114.2 |  | 6 | 19 | 3477 | 300 | 3549 |
| 114.3 |  | 38 | 3 | 237 | 3116 | 3125 |
| 114.4 |  | 114 | 1 | 229 | 26220 | 26221 |
| 115.1 | 230 | 1 | 115 | 13455 | 232 | 13457 |
| 115.2 |  | 5 | 23 | 759 | 280 | 809 |
| 115.3 |  | 23 | 5 | 255 | 1288 | 1313 |
| 115.4 |  | 115 | 1 | 231 | 26680 | 26681 |
| 116.1 | 232 | 4 | 29 | 1073 | 264 | 1105 |
| 116.2 |  | 116 | 1 | 233 | 27144 | 27145 |
| 117.1 | 234 | 1 | 117 | 13923 | 236 | 13925 |
| 117.2 |  | 9 | 13 | 403 | 396 | 565 |
| 117.3 |  | 13 | 9 | 315 | 572 | 653 |
| 117.4 |  | 117 | 1 | 235 | 27612 | 27613 |
| 118.1 | 236 | 2 | 59 | 3717 | 244 | 3725 |
| 118.2 |  | 118 | 1 | 237 | 28084 | 28085 |
| 119.1 | 238 | 1 | 119 | 14399 | 240 | 14401 |
| 119.2 |  | 7 | 17 | 527 | 336 | 625 |
| 119.3 |  | 17 | 7 | 287 | 816 | 865 |
| 119.4 |  | 119 | 1 | 239 | 28560 | 28561 |
| 120.1 | 240 | 8 | 15 | 465 | 368 | 593 |
| 120.2 |  | 24 | 5 | 265 | 1392 | 1417 |
| 120.3 |  | 40 | 3 | 249 | 3440 | 3449 |
| 120.4 |  | 120 | 1 | 241 | 29040 | 29041 |
| 121.1 | 242 | 1 | 121 | 14883 | 244 | 14885 |
| 121.2 |  | 121 | 1 | 243 | 29524 | 29525 |
| 122.1 | 244 | 2 | 61 | 3965 | 252 | 3973 |
| 122.2 |  | 122 | 1 | 245 | 30012 | 30013 |
| 123.1 | 246 | 1 | 123 | 15375 | 248 | 15377 |
| 123.2 |  | 3 | 41 | 1927 | 264 | 1945 |
| 123.3 |  | 41 | 3 | 255 | 3608 | 3617 |
| 123.4 |  | 123 | 1 | 247 | 30504 | 30505 |
| 124.1 | 248 | 4 | 31 | 1209 | 280 | 1241 |
| 124.2 |  | 124 | 1 | 249 | 31000 | 31001 |
| 125.1 | 250 | 1 | 125 | 15875 | 252 | 15877 |
| 125.2 |  | 125 | 1 | 251 | 31500 | 31501 |
| 126.1 | 252 | 2 | 63 | 4221 | 260 | 4229 |
| 126.2 |  | 14 | 9 | 333 | 644 | 725 |
| 126.3 |  | 18 | 7 | 301 | 900 | 949 |
| 126.4 |  | 126 | 1 | 253 | 32004 | 32005 |
| 127.1 | 254 | 1 | 127 | 16383 | 256 | 16385 |
| 127.2 |  | 127 | 1 | 255 | 32512 | 32513 |
| 128.1 | 256 | 128 | 1 | 257 | 33024 | 33025 |

| $N.n_i$ | $S$ | $t$ | $l$ | $x = S + l^2$ | $y = S + 2t^2$ | $a = S + l^2 + 2t^2$ |
|---|---|---|---|---|---|---|
| 129.1 | 258 | 1 | 129 | 16899 | 260 | 16901 |
| 129.2 | | 3 | **43** | 2107 | 276 | 2125 |
| 129.3 | | 43 | 3 | 267 | 3956 | 3965 |
| 129.4 | | 129 | 1 | 259 | 33540 | 33541 |
| 130.1 | 260 | 2 | 65 | 4485 | 268 | 4493 |
| 130.2 | | 10 | 13 | 429 | 460 | 629 |
| 130.3 | | 26 | 5 | 285 | 1612 | 1637 |
| 130.4 | | 130 | 1 | 261 | 34060 | 34061 |
| 131.1 | 262 | 1 | 131 | 17423 | 264 | 17425 |
| 131.2 | | 131 | 1 | 263 | 34584 | 34585 |
| 132.1 | 264 | 4 | 33 | 1353 | 296 | 1385 |
| 132.2 | | 12 | 11 | 385 | 552 | 673 |
| 132.3 | | 44 | 3 | 273 | 4136 | 4145 |
| 132.4 | | 132 | 1 | 265 | 35112 | 35113 |
| 133.1 | 266 | 1 | 133 | 17955 | 268 | 17957 |
| 133.2 | | 7 | 19 | 627 | 364 | 725 |
| 133.3 | | 19 | 7 | 315 | 988 | 1037 |
| 133.4 | | 133 | 1 | 267 | 35644 | 35645 |
| 134.1 | 268 | 2 | 67 | 4757 | 276 | 4765 |
| 134.2 | | 134 | 1 | 269 | 36180 | 36181 |
| 135.1 | 270 | 1 | 135 | 18495 | 272 | 18497 |
| 135.2 | | 5 | 27 | 999 | 320 | 1049 |
| 135.3 | | 27 | 5 | 295 | 1728 | 1753 |
| 135.4 | | 135 | 1 | 271 | 36720 | 36721 |
| 136.1 | 272 | 8 | 17 | 561 | 400 | 689 |
| 136.2 | | 136 | 1 | 273 | 37264 | 37265 |
| 137.1 | 274 | 1 | 137 | 19043 | 276 | 19045 |
| 137.2 | | 137 | 1 | 275 | 37812 | 37813 |
| 138.1 | 276 | 2 | 69 | 5037 | 284 | 5045 |
| 138.2 | | 6 | 23 | 805 | 348 | 877 |
| 138.3 | | 46 | 3 | 285 | 4508 | 4517 |
| 138.4 | | 138 | 1 | 277 | 38364 | 38365 |
| 139.1 | 278 | 1 | 139 | 19599 | 280 | 19601 |
| 139.2 | | 139 | 1 | 279 | 38920 | 38921 |
| 140.1 | 280 | 4 | 35 | 1505 | 312 | 1537 |
| 140.2 | | 20 | 7 | 329 | 1080 | 1129 |
| 140.3 | | 28 | 5 | 305 | 1848 | 1873 |
| 140.4 | | 140 | 1 | 281 | 39480 | 39481 |
| 141.1 | 282 | 1 | 141 | 20163 | 284 | 20165 |
| 141.2 | | 3 | 47 | 2491 | 300 | 2509 |
| 141.3 | | 47 | 3 | 301 | 4700 | 4709 |
| 141.4 | | 141 | 1 | 283 | 40044 | 40045 |
| 142.1 | 284 | 2 | 71 | 5325 | 292 | 5333 |
| 142.2 | | 142 | 1 | 285 | 40612 | 40613 |
| 143.1 | 286 | 1 | 143 | 20735 | 288 | 20737 |
| 143.2 | | 11 | 13 | 455 | 528 | 697 |
| 143.3 | | 13 | 11 | 407 | 624 | 745 |
| 143.4 | | 143 | 1 | 287 | 41184 | 41185 |
| 144.1 | 288 | 16 | 9 | 369 | 800 | 881 |
| 144.2 | | 144 | 1 | 289 | 41760 | 41761 |
| 145.1 | 290 | 1 | 145 | 21315 | 292 | 21317 |
| 145.2 | | 5 | 29 | 1131 | 340 | 1181 |
| 145.3 | | 29 | 5 | 315 | 1972 | 1997 |
| 145.4 | | 145 | 1 | 291 | 42340 | 42341 |
| 146.1 | 292 | 2 | 73 | 5621 | 300 | 5629 |
| 146.2 | | 146 | 1 | 293 | 42924 | 42925 |

| $N.n_i$ | $S$ | $t$ | $l$ | $x = S + l^2$ | $y = S + 2t^2$ | $a = S + l^2 + 2t^2$ |
|---|---|---|---|---|---|---|
| 147.1 | 294 | 1 | 147 | 21903 | 296 | 21905 |
| 147.2 |  | 3 | 49 | 2695 | 312 | 2713 |
| 147.3 |  | 49 | 3 | 303 | 5096 | 5105 |
| 147.4 |  | 147 | 1 | 295 | 43512 | 43513 |
| 148.1 | 296 | 4 | 37 | 1665 | 328 | 1697 |
| 148.2 |  | 148 | 1 | 297 | 44104 | 44105 |
| 149.1 | 298 | 1 | 149 | 22499 | 300 | 22501 |
| 149.2 |  | 149 | 1 | 299 | 44700 | 44701 |
| 150.1 | 300 | 2 | 75 | 5925 | 308 | 5933 |
| 150.2 |  | 6 | 25 | 925 | 372 | 997 |
| 150.3 |  | 50 | 3 | 309 | 5300 | 5309 |
| 150.4 |  | 150 | 1 | 301 | 45300 | 45301 |
| 151.1 | 302 | 1 | 151 | 23103 | 304 | 23105 |
| 151.2 |  | 151 | 1 | 303 | 45904 | 45905 |
| 152.1 | 304 | 8 | 19 | 665 | 432 | 793 |
| 152.2 |  | 152 | 1 | 305 | 46512 | 46513 |
| 153.1 | 306 | 1 | 153 | 23715 | 308 | 23717 |
| 153.2 |  | 9 | 17 | 595 | 468 | 757 |
| 153.3 |  | 17 | 9 | 387 | 884 | 965 |
| 153.4 |  | 153 | 1 | 307 | 47124 | 47125 |
| 154.1 | 308 | 2 | 77 | 6237 | 316 | 6245 |
| 154.2 |  | 14 | 11 | 429 | 700 | 821 |
| 154.3 |  | 22 | 7 | 357 | 1276 | 1325 |
| 154.4 |  | 154 | 1 | 309 | 47740 | 47741 |
| 155.1 | 310 | 1 | 155 | 24335 | 312 | 24337 |
| 155.2 |  | 5 | 31 | 1271 | 360 | 1321 |
| 155.3 |  | 31 | 5 | 335 | 2232 | 2257 |
| 155.4 |  | 155 | 1 | 311 | 48360 | 48361 |
| 156.1 | 312 | 4 | 39 | 1833 | 344 | 1865 |
| 156.2 |  | 12 | 13 | 481 | 600 | 769 |
| 156.3 |  | 52 | 3 | 321 | 5720 | 5729 |
| 156.4 |  | 156 | 1 | 313 | 48984 | 48985 |
| 157.1 | 314 | 1 | 157 | 24963 | 316 | 24965 |
| 157.2 |  | 157 | 1 | 315 | 49612 | 49613 |
| 158.1 | 316 | 2 | 79 | 6557 | 324 | 6565 |
| 158.2 |  | 158 | 1 | 317 | 50244 | 50245 |
| 159.1 | 318 | 1 | 159 | 25599 | 320 | 25601 |
| 159.2 |  | 3 | 53 | 3127 | 336 | 3145 |
| 159.3 |  | 53 | 3 | 327 | 5936 | 5945 |
| 159.4 |  | 159 | 1 | 319 | 50880 | 50881 |
| 160.1 | 320 | 32 | 5 | 345 | 2368 | 2393 |
| 160.2 |  | 160 | 1 | 321 | 51520 | 51521 |
| 161.1 | 322 | 1 | 161 | 26243 | 324 | 26245 |
| 161.2 |  | 7 | 23 | 851 | 420 | 949 |
| 161.3 |  | 23 | 7 | 371 | 1380 | 1429 |
| 161.4 |  | 161 | 1 | 323 | 52164 | 52165 |
| 162.1 | 324 | 2 | 81 | 6885 | 332 | 6893 |
| 162.2 |  | 162 | 1 | 325 | 52812 | 52813 |
| 163.1 | 326 | 1 | 163 | 26895 | 328 | 26897 |
| 163.2 |  | 163 | 1 | 327 | 53464 | 53465 |
| 164.1 | 328 | 4 | 41 | 2009 | 360 | 2041 |
| 164.2 |  | 164 | 1 | 329 | 54120 | 54121 |
| 165.1 | 330 | 1 | 165 | 27555 | 332 | 27557 |
| 165.2 |  | 3 | 55 | 3355 | 348 | 3373 |
| 165.3 |  | 5 | 33 | 1419 | 380 | 1469 |
| 165.4 |  | 11 | 15 | 555 | 572 | 797 |

| $N.n_i$ | $S$ | $t$ | $l$ | $x = S + l^2$ | $y = S + 2t^2$ | $a = S + l^2 + 2t^2$ |
|---|---|---|---|---|---|---|
| 165.5 |  | 15 | 11 | 451 | 780 | 901 |
| 165.6 |  | 33 | 5 | 355 | 2508 | 2533 |
| 165.7 |  | 55 | 3 | 339 | 6380 | 6389 |
| 165.8 |  | 165 | 1 | 331 | 54780 | 54781 |
| 166.1 | 332 | 2 | 83 | 7221 | 340 | 7229 |
| 166.2 |  | 166 | 1 | 333 | 55444 | 55445 |
| 167.1 | 334 | 1 | 167 | 28223 | 336 | 28225 |
| 167.2 |  | 167 | 1 | 335 | 56112 | 56113 |
| 168.1 | 336 | 8 | 21 | 777 | 464 | 905 |
| 168.2 |  | 24 | 7 | 385 | 1488 | 1537 |
| 168.3 |  | 56 | 3 | 345 | 6608 | 6617 |
| 168.4 |  | 168 | 1 | 337 | 56784 | 56785 |
| 169.1 | 338 | 1 | 169 | 28899 | 340 | 28901 |
| 169.2 |  | 169 | 1 | 339 | 57460 | 57461 |
| 170.1 | 340 | 2 | 85 | 7565 | 348 | 7573 |
| 170.2 |  | 10 | 17 | 629 | 540 | 829 |
| 170.3 |  | 34 | 5 | 365 | 2652 | 2677 |
| 170.4 |  | 170 | 1 | 341 | 58140 | 58141 |
| 171.1 | 342 | 1 | 171 | 29583 | 344 | 29585 |
| 171.2 |  | 9 | 19 | 703 | 504 | 865 |
| 171.3 |  | 19 | 9 | 423 | 1064 | 1145 |
| 171.4 |  | 171 | 1 | 343 | 58824 | 58825 |
| 172.1 | 344 | 4 | 43 | 2193 | 376 | 2225 |
| 172.2 |  | 172 | 1 | 345 | 59512 | 59513 |
| 173.1 | 346 | 1 | 173 | 30275 | 348 | 30277 |
| 173.2 |  | 173 | 1 | 347 | 60204 | 60205 |
| 174.1 | 348 | 2 | 87 | 7917 | 356 | 7925 |
| 174.2 |  | 6 | 29 | 1189 | 420 | 1261 |
| 174.3 |  | 58 | 3 | 357 | 7076 | 7085 |
| 174.4 |  | 174 | 1 | 349 | 60900 | 60901 |
| 175.1 | 350 | 1 | 175 | 30975 | 352 | 30977 |
| 175.2 |  | 7 | 25 | 975 | 448 | 1073 |
| 175.3 |  | 25 | 7 | 399 | 1600 | 1649 |
| 175.4 |  | 175 | 1 | 351 | 61600 | 61601 |
| 176.1 | 352 | 16 | 11 | 473 | 864 | 985 |
| 176.2 |  | 176 | 1 | 353 | 62304 | 62305 |
| 177.1 | 354 | 1 | 177 | 31683 | 356 | 31685 |
| 177.2 |  | 3 | 59 | 3835 | 372 | 3853 |
| 177.3 |  | 59 | 3 | 363 | 7316 | 7325 |
| 177.4 |  | 177 | 1 | 355 | 63012 | 63013 |
| 178.1 | 356 | 2 | 89 | 8277 | 364 | 8285 |
| 178.2 |  | 178 | 1 | 357 | 63724 | 63725 |
| 179.1 | 358 | 1 | 179 | 32399 | 360 | 32401 |
| 179.2 |  | 179 | 1 | 359 | 64440 | 64441 |
| 180.1 | 360 | 4 | 45 | 2385 | 392 | 2417 |
| 180.2 |  | 20 | 9 | 441 | 1160 | 1241 |
| 180.3 |  | 36 | 5 | 385 | 2952 | 2977 |
| 180.4 |  | 180 | 1 | 361 | 65160 | 65161 |
| 181.1 | 362 | 1 | 181 | 33123 | 364 | 33125 |
| 181.2 |  | 181 | 1 | 363 | 65884 | 65885 |
| 182.1 | 364 | 2 | 91 | 8645 | 372 | 8653 |
| 182.2 |  | 14 | 13 | 533 | 756 | 925 |
| 182.3 |  | 26 | 7 | 413 | 1716 | 1765 |
| 182.4 |  | 182 | 1 | 365 | 66612 | 66613 |
| 183.1 | 366 | 1 | 183 | 33855 | 368 | 33857 |
| 183.2 |  | 3 | 61 | 4087 | 384 | 4105 |

| $N.n_i$ | $S$ | $t$ | $l$ | $x = S + l^2$ | $y = S + 2t^2$ | $a = S + l^2 + 2t^2$ |
|---|---|---|---|---|---|---|
| 183.3 |  | 61 | 3 | 375 | 7808 | 7817 |
| 183.4 |  | 183 | 1 | 367 | 67344 | 76345 |
| 184.1 | 368 | 8 | 23 | 897 | 496 | 1025 |
| 184.2 |  | 184 | 1 | 369 | 68080 | 68081 |
| 185.1 | 370 | 1 | 185 | 34595 | 372 | 34597 |
| 185.2 |  | 5 | 37 | 1739 | 420 | 1789 |
| 185.3 |  | 37 | 5 | 395 | 3108 | 3133 |
| 185.4 |  | 185 | 1 | 371 | 68820 | 68821 |
| 186.1 | 372 | 2 | 93 | 9021 | 380 | 9029 |
| 186.2 |  | 6 | 31 | 1333 | 444 | 1405 |
| 186.3 |  | 62 | 3 | 381 | 8060 | 8069 |
| 186.4 |  | 186 | 1 | 373 | 69564 | 69565 |
| 187.1 | 374 | 1 | 187 | 35343 | 376 | 35345 |
| 187.2 |  | 11 | 17 | 663 | 616 | 905 |
| 187.3 |  | 17 | 11 | 495 | 952 | 1073 |
| 187.4 |  | 187 | 1 | 375 | 70312 | 70313 |
| 188.1 | 376 | 4 | 47 | 2585 | 408 | 2617 |
| 188.2 |  | 188 | 1 | 377 | 71064 | 71065 |
| 189.1 | 378 | 1 | 189 | 36099 | 380 | 36101 |
| 189.2 |  | 7 | 27 | 1107 | 476 | 1205 |
| 189.3 |  | 27 | 7 | 427 | 1836 | 1885 |
| 189.4 |  | 189 | 1 | 379 | 71820 | 71821 |
| 190.1 | 380 | 2 | 95 | 9405 | 388 | 9413 |
| 190.2 |  | 10 | 19 | 741 | 580 | 941 |
| 190.3 |  | 38 | 5 | 405 | 3268 | 3293 |
| 190.4 |  | 190 | 1 | 381 | 72580 | 72581 |
| 191.1 | 382 | 1 | 191 | 36863 | 384 | 36865 |
| 191.2 |  | 191 | 1 | 383 | 73344 | 73345 |
| 192.1 | 384 | 64 | 3 | 393 | 8576 | 8585 |
| 192.2 |  | 192 | 1 | 385 | 74112 | 74113 |
| 193.1 | 386 | 1 | 193 | 37635 | 388 | 37637 |
| 193.2 |  | 193 | 1 | 387 | 74884 | 74885 |
| 194.1 | 388 | 2 | 97 | 9797 | 396 | 9805 |
| 194.2 |  | 194 | 1 | 389 | 75660 | 75661 |
| 195.1 | 390 | 1 | 195 | 38415 | 392 | 38417 |
| 195.2 |  | 3 | 65 | 4615 | 408 | 4633 |
| 195.3 |  | 5 | 39 | 1911 | 440 | 1961 |
| 195.4 |  | 13 | 15 | 615 | 728 | 953 |
| 195.5 |  | 15 | 13 | 559 | 840 | 1009 |
| 195.6 |  | 39 | 5 | 415 | 3432 | 3457 |
| 195.7 |  | 65 | 3 | 399 | 8840 | 8849 |
| 195.8 |  | 195 | 1 | 391 | 76440 | 76441 |
| 196.1 | 392 | 4 | 49 | 2793 | 424 | 2825 |
| 196.2 |  | 196 | 1 | 393 | 77224 | 77225 |
| 197.1 | 394 | 1 | 197 | 39203 | 396 | 39205 |
| 197.2 |  | 197 | 1 | 395 | 78012 | 78013 |
| 198.1 | 396 | 2 | 99 | 10197 | 404 | 10205 |
| 198.2 |  | 18 | 11 | 517 | 1044 | 1165 |
| 198.3 |  | 22 | 9 | 477 | 1364 | 1445 |
| 198.4 |  | 198 | 1 | 397 | 78804 | 78805 |
| 199.1 | 398 | 1 | 199 | 39999 | 400 | 40001 |
| 199.2 |  | 199 | 1 | 399 | 79600 | 79601 |
| 200.1 | 400 | 8 | 25 | 1025 | 528 | 1153 |
| 200.2 |  | 200 | 1 | 401 | 1200 | 1201 |
| 201.1 | 402 | 1 | 201 | 40803 | 404 | 40805 |
| 201.2 |  | 3 | 67 | 4891 | 420 | 4909 |

| $N.n_i$ | $S$ | $t$ | $l$ | $x = S + l^2$ | $y = S + 2t^2$ | $a = S + l^2 + 2t^2$ |
|---|---|---|---|---|---|---|
| 201.3 |  | 67 | 3 | 411 | 9380 | 9389 |
| 201.4 |  | 201 | 1 | 403 | 81204 | 81205 |
| 202.1 | 404 | 2 | 101 | 10605 | 412 | 10613 |
| 202.2 |  | 202 | 1 | 405 | 82012 | 82013 |
| 203.1 | 406 | 1 | 203 | 41615 | 408 | 41617 |
| 203.2 |  | 7 | 29 | 1247 | 504 | 1345 |
| 203.3 |  | 29 | 7 | 455 | 2088 | 2137 |
| 203.4 |  | 203 | 1 | 407 | 82824 | 82825 |
| 204.1 | 408 | 4 | 51 | 3009 | 440 | 3041 |
| 204.2 |  | 12 | 17 | 697 | 696 | 985 |
| 204.3 |  | 68 | 3 | 417 | 9656 | 9665 |
| 204.4 |  | 204 | 1 | 409 | 83640 | 83641 |
| 205.1 | 410 | 1 | 205 | 42435 | 412 | 42437 |
| 205.2 |  | 5 | 41 | 2091 | 460 | 2141 |
| 205.3 |  | 41 | 5 | 435 | 3772 | 3797 |
| 205.4 |  | 205 | 1 | 411 | 84460 | 84461 |
| 206.1 | 412 | 2 | 103 | 11021 | 420 | 11029 |
| 206.2 |  | 206 | 1 | 413 | 85284 | 85285 |
| 207.1 | 414 | 1 | 207 | 43263 | 416 | 43265 |
| 207.2 |  | 9 | 23 | 943 | 576 | 1105 |
| 207.3 |  | 23 | 9 | 495 | 1472 | 1553 |
| 207.4 |  | 207 | 1 | 415 | 86112 | 86113 |
| 208.1 | 416 | 16 | 13 | 585 | 928 | 1097 |
| 208.2 |  | 208 | 1 | 417 | 86944 | 86945 |
| 209.1 | 418 | 1 | 209 | 44099 | 420 | 44101 |
| 209.2 |  | 11 | 19 | 779 | 660 | 1021 |
| 209.3 |  | 19 | 11 | 539 | 1140 | 1261 |
| 209.4 |  | 209 | 1 | 419 | 87780 | 87781 |
| 210.1 | 420 | 2 | 105 | 11445 | 428 | 11453 |
| 210.2 |  | 6 | 35 | 1645 | 492 | 1717 |
| 210.3 |  | 10 | 21 | 861 | 620 | 1061 |
| 210.4 |  | 14 | 15 | 645 | 812 | 1037 |
| 210.5 |  | 30 | 7 | 469 | 2220 | 2269 |
| 210.6 |  | 42 | 5 | 445 | 3948 | 3973 |
| 210.7 |  | 70 | 3 | 429 | 10220 | 10229 |
| 210.8 |  | 210 | 1 | 421 | 88620 | 88621 |
| 211.1 | 422 | 1 | 211 | 44943 | 424 | 44945 |
| 211.2 |  | 211 | 1 | 423 | 89464 | 89465 |
| 212.1 | 424 | 4 | 53 | 3233 | 456 | 3265 |
| 212.2 |  | 212 | 1 | 425 | 90312 | 90313 |
| 213.1 | 426 | 1 | 213 | 45795 | 428 | 45797 |
| 213.2 |  | 3 | 71 | 5467 | 444 | 5485 |
| 213.3 |  | 71 | 3 | 435 | 10508 | 10517 |
| 213.4 |  | 213 | 1 | 427 | 91164 | 91165 |
| 214.1 | 428 | 2 | 107 | 11877 | 436 | 11885 |
| 214.2 |  | 214 | 1 | 429 | 92020 | 92021 |
| 215.1 | 430 | 1 | 215 | 46655 | 432 | 46657 |
| 215.2 |  | 5 | 43 | 2279 | 480 | 2329 |
| 215.3 |  | 43 | 5 | 455 | 4128 | 4153 |
| 215.4 |  | 215 | 1 | 431 | 92880 | 92881 |
| 216.1 | 432 | 8 | 27 | 1161 | 560 | 1289 |
| 216.2 |  | 216 | 1 | 433 | 93744 | 93745 |
| 217.1 | 434 | 1 | 217 | 47523 | 436 | 47525 |
| 217.2 |  | 7 | 31 | 1395 | 532 | 1493 |
| 217.3 |  | 31 | 7 | 483 | 2356 | 2405 |
| 217.4 |  | 217 | 1 | 435 | 94612 | 94613 |

| $N.n_i$ | $S$ | $t$ | $l$ | $x = S + l^2$ | $y = S + 2t^2$ | $a = S + l^2 + 2t^2$ |
|---|---|---|---|---|---|---|
| 218.1 | 436 | 2 | 109 | 12317 | 444 | 12325 |
| 218.2 |  | 218 | 1 | 437 | 95484 | 95485 |
| 219.1 | 438 | 1 | 219 | 48399 | 440 | 48401 |
| 219.2 |  | 3 | 73 | 5767 | 456 | 5785 |
| 219.3 |  | 73 | 3 | 447 | 11096 | 11105 |
| 219.4 |  | 219 | 1 | 439 | 96360 | 96361 |
| 220.1 | 440 | 4 | 55 | 3465 | 472 | 3497 |
| 220.2 |  | 20 | 11 | 561 | 1240 | 1361 |
| 220.3 |  | 44 | 5 | 465 | 4312 | 4337 |
| 220.4 |  | 220 | 1 | 441 | 97240 | 97241 |
| 221.1 | 442 | 1 | 221 | 49283 | 444 | 49285 |
| 221.2 |  | 13 | 17 | 731 | 780 | 1069 |
| 221.3 |  | 17 | 13 | 611 | 1020 | 1189 |
| 221.4 |  | 221 | 1 | 443 | 98124 | 98125 |
| 222.1 | 444 | 2 | 111 | 12765 | 452 | 12773 |
| 222.2 |  | 6 | 37 | 1813 | 516 | 1885 |
| 222.3 |  | 74 | 3 | 453 | 11396 | 11405 |
| 222.4 |  | 222 | 1 | 445 | 99012 | 99013 |
| 223.1 | 446 | 1 | 223 | 50175 | 448 | 50177 |
| 223.2 |  | 223 | 1 | 447 | 99904 | 99905 |
| 224.1 | 448 | 32 | 7 | 497 | 2496 | 2545 |
| 224.2 |  | 224 | 1 | 449 | 100800 | 100801 |
| 225.1 | 450 | 1 | 225 | 51075 | 452 | 51077 |
| 225.2 |  | 9 | 25 | 1075 | 612 | 1237 |
| 225.3 |  | 25 | 9 | 531 | 1700 | 1781 |
| 225.4 |  | 225 | 1 | 451 | 101700 | 101701 |
| 226.1 | 452 | 2 | 113 | 13221 | 460 | 13229 |
| 226.2 |  | 226 | 1 | 453 | 102604 | 102605 |
| 227.1 | 454 | 1 | 227 | 51983 | 456 | 51985 |
| 227.2 |  | 227 | 1 | 455 | 103512 | 103513 |
| 228.1 | 456 | 4 | 57 | 3705 | 488 | 3737 |
| 228.2 |  | 12 | 19 | 817 | 744 | 1105 |
| 228.3 |  | 76 | 3 | 465 | 12008 | 12017 |
| 228.4 |  | 228 | 1 | 457 | 104424 | 104425 |
| 229.1 | 458 | 1 | 229 | 52899 | 460 | 52901 |
| 229.2 |  | 229 | 1 | 459 | 105340 | 105341 |
| 230.1 | 460 | 2 | 115 | 13685 | 468 | 13693 |
| 230.2 |  | 10 | 23 | 989 | 660 | 1189 |
| 230.3 |  | 46 | 5 | 485 | 4692 | 4717 |
| 230.4 |  | 230 | 1 | 461 | 106260 | 106261 |
| 231.1 | 462 | 1 | 231 | 53823 | 464 | 53825 |
| 231.2 |  | 3 | 77 | 6391 | 480 | 6409 |
| 231.3 |  | 7 | 33 | 1551 | 560 | 1649 |
| 231.4 |  | 11 | 21 | 903 | 704 | 1145 |
| 231.5 |  | 21 | 11 | 583 | 1344 | 1465 |
| 231.6 |  | 33 | 7 | 511 | 2640 | 2689 |
| 231.7 |  | 77 | 3 | 471 | 12320 | 12329 |
| 231.8 |  | 231 | 1 | 463 | 107184 | 107185 |
| 232.1 | 464 | 8 | 29 | 1305 | 592 | 1433 |
| 232.2 |  | 232 | 1 | 465 | 108112 | 108113 |
| 233.1 | 466 | 1 | 233 | 54755 | 468 | 54757 |
| 233.2 |  | 233 | 1 | 467 | 109044 | 109045 |
| 234.1 | 468 | 2 | 117 | 14157 | 476 | 14165 |
| 234.2 |  | 18 | 13 | 637 | 1116 | 1285 |
| 234.3 |  | 26 | 9 | 549 | 1820 | 1901 |
| 234.4 |  | 234 | 1 | 469 | 109980 | 109981 |

| $N.n_i$ | $S$ | $t$ | $l$ | $x = S + l^2$ | $y = S + 2t^2$ | $a = S + l^2 + 2t^2$ |
|---|---|---|---|---|---|---|
| 235.1 | 470 | 1 | 235 | 55695 | 472 | 55697 |
| 235.2 | | 5 | 47 | 2679 | 520 | 2729 |
| 235.3 | | 47 | 5 | 495 | 4888 | 4913 |
| 235.4 | | 235 | 1 | 471 | 110920 | 110921 |
| 236.1 | 472 | 4 | 59 | 3953 | 504 | 3985 |
| 236.2 | | 236 | 1 | 473 | 111864 | 111865 |
| 237.1 | 474 | 1 | 237 | 56643 | 476 | 56645 |
| 237.2 | | 3 | 79 | 6715 | 492 | 6733 |
| 237.3 | | 79 | 3 | 483 | 12956 | 12965 |
| 237.4 | | 237 | 1 | 475 | 112812 | 112813 |
| 238.1 | 476 | 2 | 119 | 14637 | 484 | 14645 |
| 238.2 | | 14 | 17 | 765 | 868 | 1157 |
| 238.3 | | 34 | 7 | 525 | 2788 | 2837 |
| 238.4 | | 238 | 1 | 477 | 113764 | 113765 |
| 239.1 | 478 | 1 | 239 | 57599 | 480 | 57601 |
| 239.2 | | 239 | 1 | 479 | 114720 | 114721 |
| 240.1 | 480 | 16 | 15 | 705 | 992 | 1217 |
| 240.2 | | 48 | 5 | 505 | 5088 | 5113 |
| 240.3 | | 80 | 3 | 489 | 13280 | 13289 |
| 240.4 | | 240 | 1 | 481 | 115680 | 115681 |
| 241.1 | 482 | 1 | 241 | 58563 | 484 | 58565 |
| 241.2 | | 241 | 1 | 483 | 116644 | 116645 |
| 242.1 | 484 | 2 | 121 | 15125 | 492 | 15133 |
| 242.2 | | 242 | 1 | 485 | 117612 | 117613 |
| 243.1 | 486 | 1 | 243 | 59535 | 488 | 59537 |
| 243.2 | | 243 | 1 | 487 | 118584 | 118585 |
| 244.1 | 488 | 4 | 61 | 4209 | 520 | 4241 |
| 244.2 | | 244 | 1 | 489 | 119560 | 119561 |
| 245.1 | 490 | 1 | 245 | 60515 | 492 | 60517 |
| 245.2 | | 5 | 49 | 2891 | 540 | 2941 |
| 245.3 | | 49 | 5 | 515 | 5292 | 5317 |
| 245.4 | | 245 | 1 | 491 | 120540 | 120541 |
| 246.1 | 492 | 2 | 123 | 15621 | 500 | 15629 |
| 246.2 | | 6 | 41 | 2173 | 564 | 2245 |
| 246.3 | | 82 | 3 | 501 | 13940 | 13949 |
| 246.4 | | 246 | 1 | 493 | 121524 | 121525 |
| 247.1 | 494 | 1 | 247 | 61503 | 496 | 61505 |
| 247.2 | | 13 | 19 | 855 | 832 | 1193 |
| 247.3 | | 19 | 13 | 663 | 1216 | 1385 |
| 247.4 | | 247 | 1 | 495 | 122512 | 122513 |
| 248.1 | 496 | 8 | 31 | 1457 | 624 | 1585 |
| 248.2 | | 248 | 1 | 497 | 123504 | 123505 |
| 249.1 | 498 | 1 | 249 | 62499 | 500 | 62501 |
| 249.2 | | 3 | 83 | 7387 | 516 | 7405 |
| 249.3 | | 83 | 3 | 507 | 14276 | 14285 |
| 249.4 | | 249 | 1 | 499 | 124500 | 124501 |
| 250.1 | 500 | 2 | 125 | 16125 | 508 | 16133 |
| 250.2 | | 250 | 1 | 501 | 125500 | 125501 |